\documentclass[11pt,a4paper,reqno,twoside]{amsart}
\usepackage[T1]{fontenc}
\usepackage{lmodern}
\usepackage{amsmath,amsthm,amssymb}
  \usepackage[margin=25mm,headheight=12pt,headsep=12pt]{geometry}
  \usepackage{booktabs,array,tabularx,longtable}
  \usepackage{cite}
  \usepackage[hidelinks,hyperfootnotes=false]{hyperref}
  \usepackage{xurl}
  \hypersetup{pdftitle={Irreducible trinomials over the binary field: an elliptic-curve criterion and explicit families},pdfauthor={Siyuan Fu, Zhiming Li, and Huaiqing Zuo}}
  \newcommand{\F}{\mathbb F}
  \newcommand{\ord}{\operatorname{ord}}
  \newcommand{\rad}{\operatorname{rad}}
  
  \theoremstyle{plain}
  \newtheorem{maintheorem}{Main Theorem}
  
  \newtheorem{theorem}{Theorem}[section]
  \newtheorem{lemma}[theorem]{Lemma}
  \newtheorem{proposition}[theorem]{Proposition}
  \newtheorem{corollary}[theorem]{Corollary}
  \theoremstyle{remark}

  \numberwithin{equation}{section}
\renewcommand{\arraystretch}{1.12}
\makeatletter
\let\arxiv@tocwrite\@tocwrite
\def\@tocwrite#1#2{\arxiv@tocwrite{#1}{#2}\def\arxiv@kind{#1}\def\arxiv@section{section}\ifx\arxiv@kind\arxiv@section\ifx\@secnumber\@empty\markboth{#2}{#2}\else\protected@edef\arxiv@heading{\sectionname\space\@secnumber.\ #2}\markboth{\arxiv@heading}{\arxiv@heading}\fi\fi}
\def\ps@headings{\let\@oddfoot\@empty\let\@evenfoot\@empty\def\@oddhead{\normalfont\scriptsize\hfil\leftmark\hfil\llap{\thepage}}\def\@evenhead{\normalfont\scriptsize\rlap{\thepage}\hfil\leftmark\hfil}\let\@mkboth\markboth}
\def\ps@references{\let\@oddfoot\@empty\let\@evenfoot\@empty\def\@oddhead{\normalfont\scriptsize\hfil References\hfil\llap{\thepage}}\def\@evenhead{\normalfont\scriptsize\rlap{\thepage}\hfil References\hfil}}
\makeatother
\title[Irreducible trinomials over the binary field]{Irreducible trinomials over the binary field: an elliptic-curve criterion and explicit families}
\author{Siyuan Fu}
\address{Weiyang College, Tsinghua University, Beijing, China}
\email{fusy24@mails.tsinghua.edu.cn}
\author{Zhiming Li}
\address{Tanwei College, Tsinghua University, Beijing, China}
\email{zm-li25@mails.tsinghua.edu.cn}
\author{Huaiqing Zuo}
\address{Department of Mathematical Sciences, Tsinghua University, Beijing 100084, China}
\email{hqzuo@mail.tsinghua.edu.cn}
\date{}

\begin{document}
\begin{abstract}
We study irreducible trinomials over $\F_2$. Using an elliptic curve
and reciprocity laws, we prove that, for $k\ge2$,
$X^{(4^k+2)/3}+X+1$ is irreducible if and only if $4^k+3$ is prime
and $2$ is a primitive root modulo that prime. By power substitution
and linearization, we also assemble explicit infinite families of
irreducible trinomials $X^m+X^n+1$ in fourteen distinct reduced
exponent classes.

\smallskip
\noindent\textit{Keywords.} Irreducible trinomial, finite field, elliptic curve,
Frobenius orbit, polynomial composition, linearized polynomial.
\par
\noindent\textit{MSC (2020).} 12E05, 11T06, 11T55, 11G20.
\end{abstract}
\maketitle

\section{Introduction}

An irreducible polynomial $f\in\F_2[X]$ of degree $m$ gives the
representation $\F_{2^m}\simeq\F_2[X]/(f)$.
Irreducible trinomials provide sparse defining polynomials for
such fields, with applications in coding theory and cryptography
\cite{LN,SunarKoc1999}.

For the fixed family $f_m(X)=X^m+X+1$,
Swan~\cite{Swan1962} gives necessary congruence restrictions on
irreducible members. Zierler~\cite{Zierler1970} lists all irreducible
members with $m\le30000$ and gives constructions from primitive
seeds. Nguyen~\cite{Nguyen2024} obtains lower bounds for the multiplicative
period of $f_m$. 
For trinomials of the more general form $X^r+X^s+1$,
Brent, Larvala, and Zimmermann~\cite{Brent2003} vary the middle
exponent $s$ at selected large prime degrees $r$ and find primitive
trinomials.
These results motivate the search for exact
irreducibility criteria on further sequences of degrees.

We study $X^m+X^n+1$, where $m>n\ge1$. We establish an
equivalence between irreducibility and a number-theoretic condition:
for $k\ge2$,
$X^{(4^k+2)/3}+X+1$ is irreducible over $\F_2$ if and only if
$4^k+3$ is prime and $2$ is a primitive root modulo that prime.
We establish this relation using an elliptic-curve description of
the roots, followed by quadratic and cubic reciprocity.
We also assemble explicit infinite families in fourteen distinct
reduced exponent classes, combining power substitution and
linearization with seed data, including nonprimitive seeds and
substitutions supported on several primes.

Using an elliptic curve, we obtain the following arithmetic
criterion for irreducibility.

\begin{maintheorem}
\label{main:arithmetic}
Let $k\ge1$, put $m=(4^k+2)/3$ and $M=4^k+3$, and set
$T_k(X)=X^m+X+1$. Then $T_k$ is irreducible over $\F_2$ if and
only if $M$ is prime and
\[
 \operatorname{lcm}\bigl(\ord_M(2),6\bigr)=M-1.
\]
\end{maintheorem}

Here the multiplicative order of $2$ modulo $M$ is defined by
\[
 \ord_M(2)=\min\{j\in\mathbb Z_{>0}:2^j\equiv1\pmod M\}.
\]
This minimum exists because $M$ is odd.
The theorem converts the polynomial irreducibility problem into a
primality condition and a multiplicative-order condition on the
integer $4^k+3$.
Main Theorem~\ref{main:arithmetic} is proved in
Subsection~\ref{sec:proof-main-a}. The proof uses the two cyclic
kernels constructed in Subsection~\ref{sec:elliptic-kernels} and
the root-orbit correspondence established in
Subsection~\ref{sec:root-orbits}.

Quadratic and cubic reciprocity give the following equivalent
formulation.

\begin{corollary}\label{cor:primitive-root}
For every integer $k\ge2$,
\begin{equation}\label{eq:primitive-root}
 \begin{gathered}
 X^{(4^k+2)/3}+X+1\text{ is irreducible over }\F_2\\
 \Longleftrightarrow\quad
 p=4^k+3\text{ is prime and }\ord_p(2)=p-1.
 \end{gathered}
\end{equation}
Thus irreducibility is equivalent to primality of $4^k+3$ together
with $2$ being a primitive root modulo that prime.
\end{corollary}

Corollary~\ref{cor:primitive-root} is proved in
Appendix~\ref{app:primitive-root}.
The restriction $k\ge2$ is necessary: for $k=1$, the polynomial
$X^2+X+1$ is irreducible, but $\ord_7(2)=3$.
Corollary~\ref{cor:primitive-root} identifies the remaining question
for this degree sequence as a number-theoretic problem about
primes and primitive roots. In particular, infinitely many of the
polynomials $T_k$ are irreducible if and only if infinitely many
integers $k\ge2$ give a prime $4^k+3$ for which $2$ is a primitive
root. Even without the primitive-root condition, it is not known
whether $4^k+3$ is prime for infinitely many integers $k$.
This is a special case of the broader problem of prime values in
exponentially growing recurrence sequences; see
Grantham and Granville~\cite{GranthamGranville2024} for conjectures
and computational evidence in this setting. The infinitude
question for $T_k$ imposes the additional requirement that $2$ be
a primitive root, and is therefore an Artin-type problem restricted
to the sparse sequence $4^k+3$.
The usual Artin primitive root conjecture, with
base $2$, concerns all primes~\cite{Moree2012}; its assertion of
infinitely many such primes does not by itself imply infinitude
under this additional restriction. The primitive root here is the
integer $2$ modulo $p$; no primitivity assertion about the
polynomial $T_k$ is made.

Elliptic-curve constructions of irreducible polynomials appear in
\cite{CL2013,BBO2024}, and related rational-map dynamics in
characteristic two are studied in \cite{Ugolini2021}.
Here we identify the roots of $T_k$ with automorphism orbits on two
cyclic kernels on $y^2+y=x^3$. The Frobenius orbit criterion
(Lemma~\ref{lem:frobenius-degree}) and the geometric torsion structure
\cite[Corollary~III.6.4(b)]{Silverman} then give Main
Theorem~\ref{main:arithmetic}; see Section~\ref{sec:elliptic}.

For the explicit constructions, Main Theorem~\ref{main:composition}
is the specialization to $\F_2$ of the classical composition
criterion~\cite{LN,genBR}, including its assertion about root orders.
We recover this criterion by a direct argument with multiplicative
orders and prime-adic valuations. It determines all admissible power
substitutions from the order of a seed polynomial.

For a monic irreducible polynomial $f\in\F_2[X]$ of degree $D\ge1$
with $f(0)=1$, choose a root $\alpha\in\F_{2^D}^{*}$ of $f$.
The multiplicative order of $\alpha$ is
\[
 e=\min\{j\in\mathbb Z_{>0}:\alpha^j=1\}.
\]
This number is independent of the chosen root and is called the
order of $f$. Equivalently, $e$ is the least positive integer $j$
such that $f(X)\mid X^j-1$.
For a positive integer $t$, let $\rad(t)$ denote the product of
its distinct prime divisors, with $\rad(1)=1$.

\begin{maintheorem}
\label{main:composition}
Let $f\in\F_2[X]$ be monic and irreducible of degree $D\ge1$,
with $f(0)=1$ and order $e$. For every positive integer $t$,
\begin{equation*}
 f(X^t)\text{ is irreducible}
 \quad\Longleftrightarrow\quad
 \rad(t)\mid e\quad\text{and}\quad
 \gcd\!\left(t,\frac{2^D-1}{e}\right)=1.
\end{equation*}
When these conditions hold, $f(X^t)$ has degree $Dt$ and order $te$.
\end{maintheorem}

The proof in Subsection~\ref{sec:composition} derives the criterion
directly over $\F_2$ without invoking the general composition theorem:
it computes the root orders and uses the lifting-the-exponent formula
to determine the root degrees. 
For an irreducible seed $X^D+X^R+1$, the criterion determines
exactly which positive integers $t$ make $X^{Dt}+X^{Rt}+1$
irreducible.

The linearization criterion associated with Ore's theory
\cite{Ore,Zierler1958,Zierler1970} gives a second construction:
a primitive trinomial $T^a+T^b+1$ yields the irreducible trinomial
$X^{2^a-1}+X^{2^b-1}+1$~\cite{Brent2003}. Proposition~\ref{thm:ore}
and Corollary~\ref{cor:ore-general} give the proof.

Table~\ref{tab:final-families} collects fourteen distinct reduced
exponent classes, each represented by an explicit infinite family
of irreducible trinomials. Proposition~\ref{thm:final-families}
determines every admissible substitution parameter in the twelve
classes A--D and gives sufficient parameter sets for the two
classes E. Each class contains infinitely many irreducible
trinomials. The classes are distinct after dividing both exponents
by their greatest common divisor and identifying reciprocal
polynomials. The two classes E arise from linearized seeds of degrees $511$
and $1023$ and require only selected prime valuations of the root
orders.
Relative to the prime-power construction recalled in
Proposition~\ref{prop:original}, the thirteen classes B--E allow
root orders divisible by at least two distinct primes. The
nonprimitive seeds also show why admissible primes are determined
by both the order and its index.

After the preliminaries in Section~\ref{sec:prelim},
Sections~\ref{sec:elliptic} and~\ref{sec:constructions} prove the
arithmetic criterion and develop the explicit constructions.
Appendix~\ref{app:certificates} records the finite seed data and
verification criteria, and Appendix~\ref{app:primitive-root} proves
Corollary~\ref{cor:primitive-root}.

\section{Notation and Frobenius orbits}\label{sec:prelim}

Unless otherwise stated, all polynomials are defined over $\F_2$. We write
\[
 f_m(X)=X^m+X+1\quad(m\ge2),\qquad
 T_{m,n}(X)=X^m+X^n+1\quad(m>n\ge1).
\]
The order of a monic polynomial $F$ with $F(0)=1$ is the least positive
integer $e$ such that $F\mid X^e-1$. For irreducible $F$ of degree $D$,
this equals the multiplicative order of each root, and $e\mid 2^D-1$.
An irreducible polynomial of degree $D$ is primitive if $e=2^D-1$.
We write $\ord_e(2)$ for the multiplicative order of $2$ modulo an odd
integer $e>1$, $\nu_p$ for the $p$-adic valuation, and $\rad(t)$ for
the product of the distinct prime divisors of $t$, with $\rad(1)=1$.

\begin{lemma}\label{lem:frobenius-degree}
The degree of an algebraic element $\alpha$ over $\F_2$ is the length of
its orbit under $\sigma:\alpha\mapsto\alpha^2$.
If $\alpha\ne0$ has multiplicative order $e>1$, this degree is
$\ord_e(2)$.
\end{lemma}
\begin{proof}
Let $r$ be the orbit length. The polynomial
$\prod_{i=0}^{r-1}(X-\alpha^{2^i})$ has coefficients fixed by squaring
and therefore belongs to $\F_2[X]$. It is divisible by the minimal
polynomial of $\alpha$, whose degree is therefore at most $r$.
For the reverse inequality, the minimal polynomial vanishes at all $r$
orbit points, so its degree is at least $r$. This proves the first assertion.
For the second, $\alpha^{2^j}=\alpha$ holds exactly when $e\mid 2^j-1$.
\end{proof}

For integers $m>n\ge1$,
put $g=\gcd(m,n)$ and define
\[
 (D,R)=(m/g,n/g),\qquad
 \mathcal C(m,n)=\{(D,R),(D,D-R)\}.
\]
Both power substitution and reciprocation preserve this class.
Consequently, infinitely many pairs $(m,n)$ in one class do not give
infinitely many distinct reduced exponent ratios.

\section{An elliptic-curve criterion for irreducibility}\label{sec:elliptic}

For an integer $k\ge1$, put
\[
 q=4^k,\qquad N=(-2)^k,\qquad M=q+3=N^2+3,
 \qquad m=\frac{q+2}{3},
\]
and write $T_k(X)=X^m+X+1$. The integer $m$ is even, so $T_k'=1$.
Moreover, $M\equiv1\pmod6$ and $\gcd(M,2N)=1$.
We prove Main Theorem~\ref{main:arithmetic} by realizing the roots
as automorphism orbits on two cyclic kernels of an elliptic curve.

\subsection{Two cyclic kernels on an elliptic curve}
\label{sec:elliptic-kernels}
Consider the elliptic curve over $\F_2$
\[
 E:\ y^2+y=x^3,
\]
with identity $O$. All point groups below consist of geometric points
in $E(\overline{\F}_2)$. Negation is given by
$-(x,y)=(x,y+1)$. Fix $\omega\in\F_4$ with
$\omega^2+\omega+1=0$, and define endomorphisms
\[
 \pi(x,y)=(x^2,y^2),\qquad \rho(x,y)=(\omega x,y),
 \qquad \phi=[1]+[2]\rho.
\]
Here $[a]$ denotes multiplication by the integer $a$, addition of
endomorphisms uses the group law on $E$, and juxtaposition denotes
composition. The chord and tangent law gives
\begin{equation}\label{el:endo}
 \rho^2+\rho+[1]=0,\qquad
 \pi^2=[-2],\qquad
 \phi^2=[-3],\qquad \pi\phi=-\phi\pi.
\end{equation}
Indeed, the points $P$, $\rho P$, and $\rho^2P$ lie on a horizontal
line and sum to $O$; the resulting identity, first obtained for a
generic point, holds as an identity of endomorphisms. The tangent at
$(x,y)$ has slope $x^2$, and the doubling formula is
\[
 [2](x,y)=(x^4,y^4+1).
\]
This proves $\pi^2=[-2]$. Expanding
$([1]+[2]\rho)^2$ and using $\rho^2+\rho+[1]=0$ gives
$\phi^2=[-3]$. Finally, $\pi\rho=\rho^2\pi$ implies
$\pi\phi=([1]+[2]\rho^2)\pi=-\phi\pi$.

For $x\ne0$, the coordinates of $\phi$ are
\begin{equation}\label{el:phi}
 \phi(x,y)=(x+x^{-2},\ y+x^{-3}+\omega).
\end{equation}
To verify this formula, first suppose that $x+\omega x^4\ne0$.
The line through $P=(x,y)$ and
$[2]\rho P=(\omega x^4,y^4+1)$ has slope
\[
 \lambda=\frac{y^4+y+1}{x+\omega x^4}
 =\frac{x^6+x^3+1}{x+\omega x^4}
 =\omega^2x^2+x^{-1}.
\]
The addition law therefore gives
\[
 \begin{aligned}
 x(\phi P)&=\lambda^2+x+\omega x^4=x+x^{-2},\\
 y(\phi P)&=\lambda\bigl(x+x(\phi P)\bigr)+y+1
           =y+x^{-3}+\omega.
 \end{aligned}
\]
The expressions on the right are regular on the open set $x\ne0$
and agree with $\phi$ on a dense open subset. Hence they give $\phi$
throughout $x\ne0$. The two points with $x=0$ are fixed by $\rho$,
so the first identity in~\eqref{el:endo} shows that they have order $3$.

Define the geometric kernel groups
\[
 K_- =\{P\in E(\overline{\F}_2):([N]-\phi)P=O\},\qquad
 K_+ =\{P\in E(\overline{\F}_2):([N]+\phi)P=O\}.
\]
For a positive integer $d$, write $E[d]=\{P:[d]P=O\}$.
We use the standard description of geometric torsion: if $d$ is odd,
then $E[d]\simeq(\mathbb Z/d\mathbb Z)^2$
\cite[Corollary~III.6.4(b), p.~86]{Silverman}.

\begin{lemma}\label{el:kernels}
The groups $K_-$ and $K_+$ are cyclic of order $M$, their intersection
is $\{O\}$, and $\pi$ interchanges them. Every nonidentity point in
either group has nonzero $x$-coordinate.
\end{lemma}
\begin{proof}
Put $a=[N]-\phi$ and $b=[N]+\phi$. By~\eqref{el:endo},
\[
 ab=ba=[M],\qquad a+b=[2N].
\]
Thus $K_-,K_+\subset E[M]$. A point in their intersection is killed by
both $[M]$ and $[2N]$, so it is $O$ because $\gcd(M,2N)=1$.
Choose an integer $c$ satisfying $2Nc\equiv1\pmod M$.
Every $P\in E[M]$ has the decomposition
\[
 P=[c]bP+[c]aP,
 \qquad [c]bP\in K_-,\quad [c]aP\in K_+.
\]
Consequently $E[M]=K_-\oplus K_+$.

The anticommutation identity in~\eqref{el:endo} gives
$\pi a=b\pi$ and $\pi b=a\pi$.
Since $\pi$ is bijective on geometric points, these identities show
that it maps each kernel bijectively onto the other. In particular,
$|K_-|=|K_+|$. The torsion description gives $|E[M]|=M^2$;
the direct sum decomposition now implies $|K_-|=|K_+|=M$.
The points with $x=0$ have order $3$ and hence cannot belong to
$E[M]$, since $3\nmid M$.

It remains to prove cyclicity. A finite abelian group is cyclic if
and only if each of its Sylow subgroups is cyclic. If a subgroup
$A\subset E[M]$ were noncyclic, the structure theorem for finite
abelian groups would therefore give a prime $\ell\mid M$ for which
the $\ell$-primary subgroup of $A$ has at least two nontrivial cyclic
factors. Its subgroup $A[\ell]$ would then have dimension at least
two over $\F_\ell$. Since $A[\ell]\subset E[\ell]$ and
$E[\ell]\simeq(\mathbb Z/\ell\mathbb Z)^2$, this forces
$A[\ell]=E[\ell]$.
If $K_-$ were noncyclic, it would thus contain $E[\ell]$ for some
$\ell\mid M$. Frobenius preserves $E[\ell]$ and is bijective on it,
so $E[\ell]=\pi(E[\ell])\subset K_+$ as well. This contradicts
$K_-\cap K_+=\{O\}$. Thus $K_-$ is cyclic, and so is its isomorphic
image $K_+$ under $\pi$.
\end{proof}

\subsection{The root-orbit correspondence}
\label{sec:root-orbits}
Let
\[
 U=\langle\rho,[-1]\rangle
 \subset\operatorname{Aut}_{\overline{\F}_2}(E,O).
\]
The commuting automorphisms $\rho$ and $[-1]$ have orders $3$ and $2$,
respectively, so $U$ has order $6$. Each element of $U$ commutes with
$\phi$ and $[N]$, and therefore preserves both kernels.
Its action on $K_\pm\setminus\{O\}$ is free. Indeed, a point fixed
by $\rho$ or $\rho^2$ has $x=0$, which is impossible by
Lemma~\ref{el:kernels}. A point fixed by $[-1]$ is killed by $[2]$,
and hence is $O$ because $M$ is odd. A point fixed by $-\rho$ or
$-\rho^2$ would also be fixed by its square, which is $\rho^2$ or
$\rho$, respectively.

For a point $P=(\beta,y)$ with $\beta\ne0$, the identity
$\pi^{2k}=[N]$ and formula~\eqref{el:phi} show that
\begin{equation}\label{el:root}
 \begin{aligned}
 P\in K_-\cup K_+
 &\ \Longleftrightarrow\ x(\pi^{2k}P)=x(\phi P)\\
 &\ \Longleftrightarrow\ \beta^q=\beta+\beta^{-2}\\
 &\ \Longleftrightarrow\ \beta^{q+2}+\beta^3+1=0.
 \end{aligned}
\end{equation}
The first equivalence uses the fact that two affine points of $E$
have the same $x$-coordinate precisely when they are equal or are
negatives of one another. Both points in that equivalence are
affine: this is immediate for $\pi^{2k}P$, and follows from
\eqref{el:phi} for $\phi P$.

Since $3m=q+2$, the map $P=(\beta,y)\mapsto\beta^3$ takes every
nonidentity point in either kernel to a root of $T_k$.
Conversely, if $\alpha$ is a root of $T_k$, then $\alpha\ne0$.
Choose $\beta\in\overline{\F}_2$ with $\beta^3=\alpha$, and choose
$y$ satisfying $y^2+y=\alpha$. Equation~\eqref{el:root} places
$(\beta,y)$ in one of the two kernels.

The fibers of this map are precisely the $U$-orbits. To see this,
suppose that $\beta'^3=\beta^3\ne0$. Then
$\beta'=\omega^j\beta$ for some $j\in\{0,1,2\}$. The corresponding
$y$-coordinates solve the same quadratic $Y^2+Y=\beta^3$, so they
are equal or differ by $1$. Thus the points lie in the same
$U$-orbit. Conversely, every element of $U$ preserves $x^3$.
Moreover, a fiber cannot meet both kernels: $U$ preserves each
kernel, and their only common point is $O$.
We have therefore obtained a bijection of sets
\begin{equation}\label{el:quotients}
 \{\alpha\in\overline{\F}_2:T_k(\alpha)=0\}
 \simeq S_-\sqcup S_+,\qquad
 S_\pm=(K_\pm\setminus\{O\})/U,
\end{equation}
where the quotients on the right denote sets of orbits.

The relations $\pi\rho=\rho^2\pi$ and
$\pi[-1]=[-1]\pi$ show that $\pi$ takes each $U$-orbit to a
$U$-orbit. Thus it induces a permutation of $S_-\sqcup S_+$, even
though $\rho$ is defined over $\F_4$. The bijection
\eqref{el:quotients} identifies this permutation with squaring on
the roots, since
\[
 x(\pi P)^3=(x(P)^3)^2.
\]
Frobenius interchanges $S_-$ and $S_+$, and its square is induced
by multiplication by $-2$ on either kernel.

\subsection{Proof of Main Theorem~\ref{main:arithmetic}}
\label{sec:proof-main-a}

\begin{proof}
Since $T_k$ is squarefree, Lemma~\ref{lem:frobenius-degree} shows
that it is irreducible precisely when Frobenius is transitive on
its roots. Under~\eqref{el:quotients}, Frobenius interchanges
$S_-$ and $S_+$. An orbit starting in $S_-$ returns to $S_-$ exactly
at its even iterates. Hence Frobenius is transitive on their union
if and only if its square is transitive on $S_-$. Equivalently,
multiplication by $-2$ must be transitive on $S_-$.

Suppose first that $M$ is composite. The cyclic group $K_-$ has
a point of order $M$ and a point of order $\ell$ for any prime
$\ell\mid M$, with $\ell<M$. Each element of $U$ preserves point
order, as does multiplication by $-2$, since $M$ is odd.
These two points therefore determine elements of $S_-$ that cannot
belong to the same orbit under multiplication by $-2$. Thus $T_k$
is reducible.

Now suppose that $M=p$ is prime. Choose a generator $G$ of $K_-$.
The map $a\mapsto[a]G$ is a bijection from $\F_p^*$ to
$K_-\setminus\{O\}$, and multiplication by an integer on the
kernel corresponds to multiplication by its residue modulo $p$.
The restriction of $\rho$ is an automorphism of the cyclic group
$K_-$, so it acts by a scalar $r\in\F_p^*$. Since $\rho^3=[1]$
and $\rho$ fixes no nonidentity point of $K_-$, the scalar $r$
has order $3$. The action of $U$ thus corresponds to multiplication
by the subgroup
\[
 H=\langle r,-1\rangle\subset\F_p^*,\qquad |H|=6.
\]
Consequently $S_-$ is identified with the coset space $\F_p^*/H$.
Multiplication by $-2$ is transitive on this space if and only if
$\langle-2,H\rangle=\F_p^*$. Since $-1\in H$, this is equivalent
to $\langle2,H\rangle=\F_p^*$. In a cyclic group, the subgroup
generated by subgroups of orders $a$ and $b$ has order
$\operatorname{lcm}(a,b)$. Applying this to $\langle2\rangle$ and
$H$ proves that transitivity is equivalent to
$\operatorname{lcm}(\operatorname{ord}_p(2),6)=p-1$, as required.
\end{proof}

For example, $k=2$ gives the irreducible polynomial $X^6+X+1$.
The reader can verify the criterion by checking
$2^{18/\ell}\not\equiv1\pmod{19}$ for $\ell=2,3$.

\section{Explicit families under power substitution}
\label{sec:constructions}

We now apply the composition and linearization criteria to
$X^m+X^n+1$ with $m>n\geq1$. Power substitution preserves the
exponent ratio, whereas linearization can change it. Starting from
primitive and nonprimitive seeds, we use their orders and indices to
determine admissible substitutions, distinguishing complete parameter
sets from sufficient ones.

\subsection{The composition criterion}\label{sec:composition}

For an irreducible polynomial of degree $D$ and nonzero constant term,
let $e$ denote the multiplicative order of any root and let
$I=(2^D-1)/e$ be its index.
For a finite set $P$ of primes, put
\[
 \mathcal S(P)=\left\{\prod_{p\in P}p^{u_p}:u_p\in\mathbb Z_{\geq0}\right\}.
\]
Thus $1\in\mathcal S(P)$, and the exponents $u_p$ vary independently.

\begin{proof}[Proof of Main Theorem~\ref{main:composition}]
For the general finite-field formulation, see \cite{LN,genBR}.
We derive the binary-field case directly from root orders and
prime-adic valuations.
The case $t=1$ is immediate. If $e=1$, then $f=X+1$, and
$X^t+1$ is reducible for $t>1$. We may therefore assume $e,t>1$.

Suppose first that the two conditions hold.  Then $t$ is odd.
For a root $\beta$ of $f(X^t)$, write $E=\ord(\beta)$.  Since
$\beta^t$ has order $e$,
\[
 E/\gcd(E,t)=e.
\]
For a prime $p\mid t$, put $a=\nu_p(E)$, $b=\nu_p(t)$, and
$c=\nu_p(e)>0$.  The preceding identity gives $a-\min(a,b)=c$,
and hence $a=b+c$. For $p\nmid t$, it gives $\nu_p(E)=\nu_p(e)$.
Thus $E=te$. Since $\ord_e(2)=D$, every positive integer $j$
for which $te\mid2^j-1$ has the form $j=Du$.  For each prime $p\mid t$,
the lifting-the-exponent formula gives
\[
 \nu_p(2^{Du}-1)=\nu_p(2^D-1)+\nu_p(u)
               =\nu_p(e)+\nu_p(u).
\]
Here $p$ is odd, and the last equality uses $p\nmid I$.
Consequently $te\mid2^{Du}-1$ if and only if $t\mid u$.
Thus $\ord_{te}(2)=Dt$, which is the degree of $\beta$ over
$\F_2$, proving irreducibility and the asserted order.

Conversely, suppose that $f(X^t)$ is irreducible.  If $t$ were even,
then $f(X^t)=f(X^{t/2})^2$, so $t$ must be odd.
Let $p$ be a prime divisor of $t$, and let $\alpha$ be a root of $f$
in $K=\F_{2^D}$.
If $p\nmid e$, choose $v$ with $pv\equiv1\pmod e$; then
$(\alpha^v)^p=\alpha$.  If instead $p\mid e$ and $p\mid I$,
then $e\mid(2^D-1)/p$.  The subgroup of $p$th powers in the cyclic
group $K^*$ has order $(2^D-1)/p$, and hence contains $\alpha$.
In either case, $f(X^p)$ has a root in $K$ and hence a factor of
degree at most $D<Dp$.  Substituting $X^{t/p}$ in the resulting
nontrivial factorization would make $f(X^t)$ reducible.
Both cases contradict irreducibility, proving necessity.
\end{proof}

\begin{corollary}
\label{gen:normalization}
For integers $m>n\ge1$, put $g=\gcd(m,n)$, $D=m/g$, and $R=n/g$.
Then $X^m+X^n+1$ is irreducible if and only if
$X^D+X^R+1$ is irreducible and its order $e$ satisfies
$\rad(g)\mid e$ and $\gcd(g,(2^D-1)/e)=1$.
\end{corollary}

\begin{proof}
A nontrivial factorization of $X^D+X^R+1$ remains nontrivial after
substitution of $X^g$. When the seed is irreducible, the assertion
follows from Main Theorem~\ref{main:composition}.
\end{proof}

\begin{proposition}\label{prop:original}
Let $p$ be an odd prime and let $h\ne X+1$ be a monic irreducible
factor of $X^p-1$ over $\F_2$.  If
$\ord_{p^2}(2)\ne\ord_p(2)$, then $h(X^{p^k})$ is irreducible
for every integer $k\geq0$, with order $p^{k+1}$ and degree
$p^k\ord_p(2)$.
\end{proposition}

\begin{proof}
Put $D=\ord_p(2)$.  The order of $h$ is $p$, and the hypothesis
is equivalent to $\nu_p(2^D-1)=1$.  Thus
$p\nmid(2^D-1)/p$, and Main Theorem~\ref{main:composition} applies
with $t=p^k$.
\end{proof}

Taking $p=3$ and $h=X^2+X+1$ gives
$X^{2\cdot3^k}+X^{3^k}+1$ for $k\geq0$.
All polynomials in Proposition~\ref{prop:original} have prime-power
order, even when the seeds are nonprimitive. We next obtain families
whose orders are divisible by several distinct primes.

\subsection{Primitive and nonprimitive seeds}

We write the families in Table~\ref{tab:final-families} in the form
\[
 (m,n)=(Dt,Rt),\qquad t\in\mathcal S(P).
\]
The labels distinguish the seed data and the construction: A is the
prime-power baseline; B uses primitive seeds of composite order;
C uses nonprimitive seeds with $\gcd(e,I)=1$; D uses a seed with
$\gcd(e,I)>1$; and E combines linearization with power substitution.
The primes in the displayed product for $t$ form the set $P$.

Table~\ref{tab:seed-orders} in Appendix~\ref{app:certificates}
(Subsection~\ref{app:seed-data}) gives the exact seed orders and
indices. The B seeds are primitive, so $I=1$. The C seeds are
nonprimitive but satisfy $\gcd(e,I)=1$. In both cases, Main
Theorem~\ref{main:composition} permits every prime divisor of $e$
in $t$, with arbitrary nonnegative exponent.

The reader can verify these finite data using
Lemma~\ref{lem:certificate} in Appendix~\ref{app:certificates}:
compute the powers in \eqref{eq:rabin} and \eqref{eq:exact-order}
by repeated squaring modulo the seed, and evaluate the required
polynomial gcds.

For D$_1$, the seed $H=X^{60}+X^{17}+1$ has index $I=1155$ and
$\gcd(e,I)=15$; see Table~\ref{tab:seed-orders}.
Thus $3$ and $5$ are excluded from $t$, while the other seven
prime divisors of $e$ are admissible by Main
Theorem~\ref{main:composition}. This gives precisely the parameter
set in row D$_1$ of Table~\ref{tab:final-families}.

\subsection{Classical linearized lifting}\label{sec:lift}

Power substitution preserves the ratio of the positive exponents.
To obtain further ratios, we use the classical construction due to Ore \cite[Chapter~III, Theorem~1, p.~260]{Ore}; see also
\cite{Zierler1958,Zierler1970}. Its trinomial specialization appears in
\cite[p.~1450]{Brent2003}. We give a direct proof of the equivalence used below.

For a monic polynomial $h(T)=\sum_{i=0}^a h_iT^i\in\F_2[T]$ of
degree $a\ge1$ with $h_0=1$, define
\[
 L_h(X)=\sum_{i=0}^a h_iX^{2^i},\qquad B_h(X)=L_h(X)/X.
\]
Thus $\deg B_h=2^a-1$. Since all coefficients lie in $\F_2$,
expansion gives the identity
\begin{equation}\label{eq:linearized-composition}
 L_{uv}=L_u\circ L_v.
\end{equation}

\begin{proposition}\label{thm:ore}
For every monic polynomial $h\in\F_2[T]$ of positive degree with
$h(0)=1$,
\begin{equation}\label{eq:ore-equivalence}
 B_h\text{ is irreducible}\quad\Longleftrightarrow\quad
 h\text{ is primitive}.
\end{equation}
\end{proposition}

\begin{proof}
Write $a=\deg h$. Additivity and $L_h'=1$ imply that the $2^a$
distinct roots of $L_h$ form an $a$-dimensional $\F_2$-vector space
$V$ in $\overline{\F}_2$. The identity $L_h(v^2)=L_h(v)^2$ shows
that squaring induces an invertible linear map $\sigma$ on $V$.
For every $v\in V$,
\[
 h(\sigma)v=\sum_{i=0}^a h_i v^{2^i}=L_h(v)=0.
\]

Suppose first that $h$ is irreducible, and let $e$ be its order.
For $0\ne v\in V$, let $c_v$ be the monic polynomial of least degree
such that $c_v(\sigma)v=0$. Division with remainder shows that $c_v$
divides every polynomial annihilating $v$, and in particular divides
$h$. Since $c_v\ne1$ and $h$ is irreducible, $c_v=h$. It follows that
\[
 \sigma^j(v)=v\quad\Longleftrightarrow\quad h\mid T^j-1
 \quad\Longleftrightarrow\quad e\mid j.
\]
Thus the minimal polynomial of any nonzero root of $L_h$ has degree
$e$, by Lemma~\ref{lem:frobenius-degree}. Such a minimal polynomial
divides $B_h$, whose degree is $2^a-1$. Consequently, for irreducible
$h$, the polynomial $B_h$ is irreducible if and only if
$e=2^a-1$, which is exactly the primitivity of $h$.

It remains to show that irreducibility of $B_h$ forces irreducibility
of $h$. If $h=uv$ with both factors of positive degree, then
$u(0)=v(0)=1$, and \eqref{eq:linearized-composition} gives
\[
 B_h(X)=B_v(X)\,B_u(L_v(X)).
\]
For $b=\deg v$, the two factors have positive degrees $2^b-1$ and
$2^a-2^b$. Thus $B_h$ is reducible whenever $h$ is reducible.
Combining this observation with the preceding argument proves both
directions of \eqref{eq:ore-equivalence}.
\end{proof}

\subsection{Linearization followed by power substitution}

We now use Proposition~\ref{thm:ore} to change the exponent ratio
before applying power substitution.

\begin{corollary}\label{cor:ore-general}
If $T^a+T^b+1$ is primitive over $\F_2$, where $a>b\geq1$,
then $X^{2^a-1}+X^{2^b-1}+1$ is irreducible.
The greatest common divisor of its positive exponents is
$2^{\gcd(a,b)}-1$.
\end{corollary}

\begin{proof}
Irreducibility follows from Proposition~\ref{thm:ore}. The identity
$\gcd(2^a-1,2^b-1)=2^{\gcd(a,b)}-1$
follows from the Euclidean algorithm.
\end{proof}

The resulting seed need not be primitive. The following criterion
determines admissible substitution primes without its full order.

\begin{corollary}
\label{cor:partial-order}\label{thm:partial-order}
Let $F\in\F_2[X]$ be irreducible of degree $D\geq1$, with
$F(0)=1$, and put $Q=2^D-1$.  For every positive integer $t$,
$F(X^t)$ is irreducible if and only if every prime divisor $p$
of $t$ satisfies
\[
 p\mid Q,\qquad X^{Q/p}\not\equiv1\pmod F.
\]
In particular, if a finite set $P$ of primes satisfies these two
conditions, then $F(X^t)$ is irreducible for every
$t\in\mathcal S(P)$.
\end{corollary}

\begin{proof}
Let $e$ be the order of $F$.  For a prime $p\mid Q$, the
divisibility $e\mid Q$ gives
\[
 X^{Q/p}\not\equiv1\pmod F
 \quad\Longleftrightarrow\quad e\nmid Q/p
 \quad\Longleftrightarrow\quad\nu_p(e)=\nu_p(Q).
\]
This is equivalent to $p\mid e$ and $p\nmid Q/e$.
If $p\nmid Q$, then $p\nmid e$, so the composition criterion fails.
The equivalence now follows from Main Theorem~\ref{main:composition};
for $t=1$, the tests are vacuous.
\end{proof}

For a specified $t$, the only factorization data required by this
criterion are its prime divisors; a complete factorization
of $Q$ and the full root order are unnecessary. A set of primes
satisfying the criterion gives a sufficient family. Completeness
requires excluding all remaining primes.

\begin{corollary}\label{cor:two-stage}
The following trinomials are irreducible for all nonnegative
integers $u,v,w,z$:
\[
\begin{aligned}
 &X^{511t}+X^{15t}+1,&\quad t&=127^u439^v,\\
 &X^{1023t}+X^{7t}+1,&\quad t&=7^u23^v89^w2147483647^z.
\end{aligned}
\]
These parameter sets are sufficient; their completeness is not asserted.
\end{corollary}

\begin{proof}
Apply Corollary~\ref{cor:ore-general} to the primitive seeds B$_3$
and B$_4$ in Table~\ref{tab:seed-orders} of
Appendix~\ref{app:certificates}. The required divisibility and
noncongruence checks are recorded in Appendix~\ref{app:certificates},
Subsection~\ref{app:two-stage}, and can be verified by repeated
squaring modulo the resulting seeds. Corollary~\ref{cor:partial-order}
then applies.
\end{proof}

\subsection{Summary and distinction of the families}
\label{sec:final-families}

Recall the reduced exponent class defined in Section~\ref{sec:prelim}:
\[
 \mathcal C(m,n)=\{(D,R),(D,D-R)\},\qquad
 D=m/g,\ R=n/g,\ g=\gcd(m,n).
\]
Both common power substitution and reciprocation preserve this class.
Table~\ref{tab:final-families} includes the reciprocal forms.
In each product for $t$, the exponents vary independently over the
nonnegative integers. In particular, $t=1$ is permitted, while the
parameter $r$ is restricted to the indicated set.

\begingroup
\small
\renewcommand{\arraystretch}{1.25}
\setlength{\tabcolsep}{4pt}
\setlength{\LTpre}{\medskipamount}
\setlength{\LTpost}{\medskipamount}
\begin{longtable}{@{}cll@{}}
\caption{Explicit irreducible trinomials and their reciprocals.
The conditions on $t$ are necessary and sufficient in A--D,
and sufficient in E.}\label{tab:final-families}\\
\toprule
Family & Polynomial and middle exponent & Parameter $t$\\
\midrule
\endfirsthead
\multicolumn{3}{c}{\tablename\ \thetable\ (continued)}\\
\toprule
Family & Polynomial and middle exponent & Parameter $t$\\
\midrule
\endhead
\midrule
\endfoot
\bottomrule
\endlastfoot
$A_1$ & $X^{2t}+X^t+1$ & $3^u$\\
$B_1$ & $X^{4t}+X^{rt}+1,\ r\in\{1,3\}$ & $3^u5^v$\\
$B_2$ & $X^{6t}+X^{rt}+1,\ r\in\{1,5\}$ & $3^u7^v$\\
$B_3$ & $X^{9t}+X^{rt}+1,\ r\in\{4,5\}$ & $7^u73^v$\\
$B_4$ & $X^{10t}+X^{rt}+1,\ r\in\{3,7\}$ & $3^u11^v31^w$\\
$B_5$ & $X^{11t}+X^{rt}+1,\ r\in\{2,9\}$ & $23^u89^v$\\
$B_6$ & $X^{15t}+X^{rt}+1,\ r\in\{1,14\}$ & $7^u31^v151^w$\\
$B_7$ & $X^{20t}+X^{rt}+1,\ r\in\{3,17\}$ & $3^a5^b11^c31^d41^h$\\
$C_1$ & $X^{12t}+X^{rt}+1,\ r\in\{5,7\}$ & $3^u7^v13^w$\\
$C_2$ & $X^{14t}+X^{rt}+1,\ r\in\{5,9\}$ & $43^u127^v$\\
$C_3$ & $X^{28t}+X^{rt}+1,\ r\in\{1,27\}$ & $29^u43^v113^w127^z$\\
$D_1$ & $X^{60t}+X^{rt}+1,\ r\in\{17,43\}$ & $\begin{gathered}13^a31^b41^c61^d\\[-0.3em]{}\cdot151^h331^i1321^j\end{gathered}$\\
$E_1$ & $X^{511t}+X^{rt}+1,\ r\in\{15,496\}$ & $127^u439^v$\\
$E_2$ & $X^{1023t}+X^{rt}+1,\ r\in\{7,1016\}$ & $7^u23^v89^w(2147483647)^z$\\
\end{longtable}
\endgroup

The following proposition summarizes the table and compares the
resulting families with the prime-power baseline.

\begin{proposition}\label{thm:final-families}
For each of the twelve rows A--D of
Table~\ref{tab:final-families} and each positive integer $t$,
each displayed trinomial is irreducible if and only if $t$ has
the product form specified in that row.
The two rows E give sufficient parameter sets.
All fourteen rows define distinct classes $\mathcal C(m,n)$,
and every row contains infinitely many irreducible trinomials.
The thirteen rows B--E lie outside the output range of
Proposition~\ref{prop:original}.
\end{proposition}

\begin{proof}
The exact orders and indices in Table~\ref{tab:seed-orders}
of Appendix~\ref{app:certificates}, together with Main
Theorem~\ref{main:composition}, determine all admissible primes
for A--D, proving both directions for these twelve rows. For E,
the assertion follows from Corollary~\ref{cor:two-stage}.

In every row, $\gcd(D,r)=1$ for each listed $r$ (with $r=1$ in A$_1$).
The fourteen reduced degrees $D$ are distinct and invariant under
power substitution and reciprocation, so the classes are distinct.
The reciprocal of $F=X^{Dt}+X^{Rt}+1$ is
$X^{Dt}F(X^{-1})=X^{Dt}+X^{(D-R)t}+1$.
Taking reciprocals of a nontrivial factorization of either polynomial
gives a nontrivial factorization of the other, with degrees preserved
because the constant terms are nonzero.
Root orders, and hence admissible substitution parameters, are also
preserved. An unrestricted prime exponent in each row gives infinitely
many distinct irreducible trinomials.

Every B--D seed has order divisible by at least two distinct primes;
the noncongruences used in Corollary~\ref{cor:two-stage}
give the same property for E.
Admissible substitution multiplies the order by $t$, preserving these
divisors. Thus these families lie outside Proposition~\ref{prop:original},
which gives only prime-power orders. Only A$_1$ is self-reciprocal,
giving $27$ displayed forms.
\end{proof}

\section{Conclusion}

Corollary~\ref{cor:primitive-root}, proved from Main
Theorem~\ref{main:arithmetic} in Appendix~\ref{app:primitive-root},
identifies the
irreducibility of $X^{(4^k+2)/3}+X+1$ for $k\ge2$ with the
primality of $4^k+3$ and the requirement that $2$ be a primitive
root modulo that prime. This equivalence permits further progress
through number-theoretic analysis. Infinitude remains open, even
for the prime values of $4^k+3$ without the primitive-root condition.

We also give a direct derivation of the classical composition
criterion over $\F_2$ through root orders and prime-adic valuations.
Together with linearization, this criterion yields the explicit
infinite families collected in Table~\ref{tab:final-families}, which
represent fourteen distinct reduced exponent classes. The parameter
sets are necessary and sufficient in A--D and sufficient in E, as
proved in Proposition~\ref{thm:final-families}. The finite seed data
and verification criteria are recorded in
Appendix~\ref{app:certificates}.

\appendix
\renewcommand{\thetheorem}{\Alph{section}.\arabic{theorem}}
\renewcommand{\theequation}{\Alph{section}.\arabic{equation}}
\section{Finite algebraic certificates}\label{app:certificates}

This appendix records the finite seed data used in
Section~\ref{sec:constructions}. The authors have checked all finite
computations in this appendix using exact arithmetic and the criteria
in Lemma~\ref{lem:certificate} and Corollary~\ref{cor:partial-order}.

\subsection{Irreducibility and order tests}\label{app:tests}

We use Rabin's irreducibility criterion~\cite{Rabin}.
The accompanying order test is the standard multiplicative-order
criterion \cite{LN}.

\begin{lemma}\label{lem:certificate}
Let $F\in\F_2[X]$ be monic of degree $D\ge1$, with $F(0)=1$.
Then $F$ is irreducible if and only if
\begin{equation}\label{eq:rabin}
 X^{2^D}\equiv X\pmod F,\qquad
 \gcd\bigl(F,X^{2^{D/\ell}}-X\bigr)=1
 \quad\text{for every prime }\ell\mid D.
\end{equation}
If $F$ is irreducible and $e$ is a positive divisor of $2^D-1$,
then $F$ has order $e$ if and only if
\begin{equation}\label{eq:exact-order}
 X^e\equiv1\pmod F,\qquad
 X^{e/p}\not\equiv1\pmod F
 \quad\text{for every prime }p\mid e.
\end{equation}
\end{lemma}

\begin{proof}
The polynomial $X^{2^D}-X$ is squarefree, since its derivative is $1$.
By Lemma~\ref{lem:frobenius-degree}, its irreducible factors are
precisely the monic irreducible polynomials whose degrees divide $D$.
The first congruence in \eqref{eq:rabin} therefore makes $F$
squarefree, with every irreducible factor having degree dividing $D$.
Every proper divisor of $D$ divides $D/\ell$ for some prime
$\ell\mid D$. The gcd conditions exclude all such degrees, so $F$
is irreducible. The converse follows from the same description of
the factors of $X^{2^D}-X$.

Let $e_0$ be the order of $F$. The first congruence in
\eqref{eq:exact-order} says $e_0\mid e$. If $e_0<e$, then
$e_0\mid e/p$ for a prime divisor $p$ of $e/e_0$, contrary to the
remaining conditions. Conversely, an element of order $e$ satisfies
all the conditions in \eqref{eq:exact-order}.
\end{proof}

\subsection{Seed orders for the substitution families}\label{app:seed-data}

For each pair $(D,R)$ in Table~\ref{tab:seed-orders}, we used
\eqref{eq:rabin} to verify that $F=X^D+X^R+1$ is irreducible.
We independently computed its order $e$ and verified it using
\eqref{eq:exact-order}. The index is $I=(2^D-1)/e$, and the seeds
with $I=1$ are primitive.
Our results agree with those of
Zierler and Brillhart~\cite{ZierlerBrillhart1968}. Their Table~1
identifies the eight seeds in our groups A and B as primitive, and
their Table~2 (p.~549) gives the periods and indices of the four
nonprimitive seeds in our groups C and D.

\begin{table}[htbp]
\centering
\small
\caption{Exact orders and indices of the twelve seeds.}
\label{tab:seed-orders}
\begin{tabular}{@{}crrrr@{}}
\toprule
Family & $D$ & $R$ & Order $e$ & Index $I$\\
\midrule
$A_1$ & 2 & 1 & 3 & 1\\
$B_1$ & 4 & 1 & 15 & 1\\
$B_2$ & 6 & 1 & 63 & 1\\
$B_3$ & 9 & 4 & 511 & 1\\
$B_4$ & 10 & 3 & 1023 & 1\\
$B_5$ & 11 & 2 & 2047 & 1\\
$B_6$ & 15 & 1 & 32767 & 1\\
$B_7$ & 20 & 3 & 1048575 & 1\\
$C_1$ & 12 & 5 & 819 & 5\\
$C_2$ & 14 & 5 & 5461 & 3\\
$C_3$ & 28 & 1 & 17895697 & 15\\
$D_1$ & 60 & 17 & 998200436889045 & 1155\\
\bottomrule
\end{tabular}
\end{table}

\subsection{The two families obtained by linearization}\label{app:two-stage}

Put $F_{511}=X^{511}+X^{15}+1$ and
$F_{1023}=X^{1023}+X^7+1$.
Their irreducibility follows from Corollary~\ref{cor:ore-general}
and the primitive seeds $(D,R)=(9,4)$ and $(10,3)$ in
Table~\ref{tab:seed-orders}.

For each pair $(D,p)$ below, $p$ is prime and divides $Q_D=2^D-1$.
Let $R_{D,p}$ be the remainder of $X^{Q_D/p}$ modulo $F_D$.
The last two columns specify an exponent $j>0$ such that
$[X^j]R_{D,p}=1$,
where $[X^j]R_{D,p}$ denotes the coefficient of $X^j$.
Thus every displayed remainder differs from $1$, as required by
Corollary~\ref{cor:partial-order}.

\begin{center}
\small
\begin{tabular}{@{}rrrr@{}}
\toprule
$D$ & $p$ & $j$ & $[X^j]R_{D,p}$\\
\midrule
511 & 127 & 1 & 1\\
511 & 439 & 1 & 1\\
1023 & 7 & 1 & 1\\
1023 & 23 & 5 & 1\\
1023 & 89 & 2 & 1\\
1023 & 2147483647 & 1 & 1\\
\bottomrule
\end{tabular}
\end{center}

Readers may verify primality by trial division, divisibility by
integer modular exponentiation, and the coefficients by repeated
squaring modulo $F_D$.
Only these six prime divisors are tested; no complete root order
is asserted for either seed.

\section{The primitive-root formulation}\label{app:primitive-root}

This appendix proves Corollary~\ref{cor:primitive-root}. We first
show that $2$ is a cubic nonresidue modulo every prime of the form
$4^k+3$. Throughout this appendix, $\zeta$ denotes a primitive
cube root of unity in characteristic zero, and
$\mathcal O=\mathbb Z[\zeta]$ is the ring of Eisenstein integers.
Its norm is
$\operatorname{Nm}(a+b\zeta)=a^2-ab+b^2$.

We use the following normalization of cubic reciprocity
\cite[Chapter~9, \S3, Theorem~1, p.~114]{IrelandRosen1990}.
A prime element $\theta$ of $\mathcal O$ is called primary if
$\theta\equiv2\pmod{3\mathcal O}$. For $\theta\nmid3$ and
$a\in\mathcal O$ with $a\not\equiv0\pmod\theta$, the cubic residue symbol
$\left(\frac a\theta\right)_3\in\{1,\zeta,\zeta^2\}$ is
characterized by
\[
 a^{(\operatorname{Nm}(\theta)-1)/3}
 \equiv\left(\frac a\theta\right)_3\pmod\theta.
\]
For primary prime elements $\theta_1,\theta_2$ of distinct norms,
neither dividing $3$, cubic reciprocity states that
\[
 \left(\frac{\theta_1}{\theta_2}\right)_3
 =
 \left(\frac{\theta_2}{\theta_1}\right)_3.
\]

\begin{lemma}\label{lem:cubic-nonresidue}
Let $k\ge1$ and suppose that $p=4^k+3$ is prime. Then
\[
 2^{(p-1)/3}\not\equiv1\pmod p.
\]
Consequently, if $d=\ord_p(2)$, then
$\nu_3(d)=\nu_3(p-1)\ge1$.
\end{lemma}

\begin{proof}
Put $N=(-2)^k$, so that $N\equiv1\pmod3$, and set
\[
 \Pi=(N+1)+(N-1)\zeta\in\mathcal O.
\]
Its norm is
\[
 \operatorname{Nm}(\Pi)
 =(N+1)^2-(N+1)(N-1)+(N-1)^2
 =N^2+3=p.
\]
Thus $\Pi$ is a prime element and
$\mathcal O/(\Pi)\simeq\F_p$. Moreover,
$\Pi\equiv2\pmod{3\mathcal O}$, so $\Pi$ is primary.
The element $2$ is also a primary prime element: indeed,
\[
 \mathcal O/(2)\simeq\F_2[Z]/(Z^2+Z+1)\simeq\F_4.
\]
Its norm is $4$, which differs from $p\ge7$, and neither $2$ nor
$\Pi$ divides $3$. Since $N$ is even,
\[
 \Pi\equiv1+\zeta\equiv\zeta^2\pmod{2\mathcal O}.
\]
The exponent defining the cubic residue symbol modulo $2$ is
$(4-1)/3=1$. Therefore
$\left(\frac{\Pi}{2}\right)_3=\zeta^2$, and cubic reciprocity gives
\begin{equation}\label{eq:cubic-character}
 2^{(p-1)/3}
 \equiv\left(\frac2\Pi\right)_3
 =\left(\frac\Pi2\right)_3
 =\zeta^2\pmod\Pi.
\end{equation}
The image of $\zeta$ in $\mathcal O/(\Pi)$ has order $3$:
it satisfies $Z^2+Z+1=0$, and cannot equal $1$ since $p\ne3$.
Hence the right-hand side of~\eqref{eq:cubic-character} is not
$1$. Under the canonical identification of the residue field with
$\F_p$, this proves the asserted noncongruence.

Finally, $d\mid p-1$. If $\nu_3(d)<\nu_3(p-1)$, then
$d\mid(p-1)/3$, contrary to the noncongruence. Since
$p\equiv1\pmod3$, this gives
$\nu_3(d)=\nu_3(p-1)\ge1$.
\end{proof}

\begin{proof}[Proof of Corollary~\ref{cor:primitive-root}]
By Main Theorem~\ref{main:arithmetic}, it suffices to consider a
prime $p=4^k+3$. Write $d=\ord_p(2)$. For $k\ge2$ we have
$p\equiv3\pmod8$. The supplementary law of quadratic reciprocity
and Euler's criterion give
\[
 2^{(p-1)/2}
 \equiv\left(\frac2p\right)
 =(-1)^{(p^2-1)/8}=-1\pmod p,
\]
where $\left(\frac2p\right)$ is the Legendre symbol.
Since $d\mid p-1$ and $\nu_2(p-1)=1$, this implies $2\mid d$.
Lemma~\ref{lem:cubic-nonresidue} gives $3\mid d$. Thus $6\mid d$
and
\[
 \operatorname{lcm}(d,6)=d.
\]
The condition in Main Theorem~\ref{main:arithmetic} is consequently
equivalent to $d=p-1$, which means that $2$ is a primitive root
modulo $p$.

\end{proof}

\section*{Acknowledgement}
\markboth{Acknowledgement}{Acknowledgement}
ChatGPT and Codex assisted with literature searches, proof checks,
verification code, and language editing. The authors take
responsibility for all content.
\par

The authors thank Professor Hao Chen for suggesting this problem.
H.~Zuo was supported by NSFC Grant 12671056 and BJNSF Grant 1252009.

\pagestyle{references}
{\footnotesize
\urlstyle{same}

}
\end{document}